\documentclass[11pt]{article}
\usepackage{amssymb}
\usepackage{amsmath}
\newtheorem{theo}{Theorem}[section]

\newtheorem{prop}[theo]{Proposition}
\newtheorem{lemm}[theo]{Lemma}
\newtheorem{coro}[theo]{Corollary}
\newtheorem{rema}[theo]{Remark}
\newtheorem{Defi}[theo]{Definition}

\newtheorem{conj}[theo]{Conjecture}

\newcommand{\cqfd}
{%
\mbox{}%
\nolinebreak%
\hfill%
\rule{2mm}{2mm}%
\medbreak%
\par%
}
\newfont{\gothic}{eufb10}
\date{\empty}
\begin{document}
\title{Rationally connected $3$-folds and symplectic geometry}
\author{Claire Voisin\\
CNRS and IH\'{E}S} \maketitle \setcounter{section}{-1}
\begin{flushright}{\it Pour Jean-Pierre Bourguignon, \`a   l'occasion de ses 60 ans}
\end{flushright}
\section{Introduction}

Let $X$ be a compact K\"ahler manifold. Denoting by $J$ the operator of
 complex structure acting on $T_X$,  K\"ahler forms on $X$ are
symplectic forms which satisfy the compatibility conditions
$$\omega(Ju,Jv)=\omega(u,v),\,u,\,v\in T_{X,x},\,\omega(u,Ju)>0,\,0\not=u\in T_{X,x}.$$
The first condition tells that $\omega$ is of type $(1,1)$.
The last condition is called the taming condition. The set of K\"ahler forms
is a convex cone, in particular connected, and thus determines a deformation class
of symplectic forms on $X$.

Let $X$ and $Y$ be two complex projective or compact K\"ahler manifolds.
We will say that $X$ and $Y$ are symplectically equivalent if for some
 symplectic  forms $\alpha$ on $X$, resp. $\beta$ on $Y$, which are
 in the deformation class of a K\"ahler form on $X$, resp. $Y$,
there is a
diffeomorphism
$$\psi: X\cong Y,$$
such that $\psi^*\beta=\alpha$. Notice that $\psi^*$ induces a bijection between the
sets of symplectic forms which are in the deformation class
of a symplectic form on $Y$ and $X$, and thus we may assume that $\alpha$ is a taming
form, or even a K\"ahler form on $X$.

In the sequel, the compact K\"ahler manifolds $X$ we will consider
 are {\it uniruled} manifolds, which means the following (cf \cite{kollarbook}):
\begin{Defi} A projective complex manifold (or compact K\"ahler) is uniruled if there exist
 compact complex manifolds $Z$ and $B$, and dominating morphisms
 $$f:Z\rightarrow X,\,g:Z\rightarrow B,$$
 where $f$ is non constant on the fibers of $g$ and  the generic fiber of $g$ is isomorphic to
 $\mathbb{P}^1$.
 \end{Defi}
 In other words, there is a (maybe singular) rational curve in $X$ passing through any point
 of $X$, where a (singular) rational curve is defined as a connected curve whose normalization
has only rational components.

 The starting point of this work is the following result, due independently
  to Koll\'ar  \cite{kollarvrai} and Ruan  \cite{ruan1}
 (we refer to \cite{ruan2}, \cite{mcduff1}, \cite{mcduff2} for purely symplectic characterizations
 and studies  of
 uniruledness) :
 \begin{theo}\label{thkr} Let $X$ and $Y$ be two symplectically equivalent compact K\"ahler
 manifolds.
 Then if $X$ is uniruled, $Y$ is also uniruled.
\end{theo}
We sketch later on the proof of this result, in order to point out why the proof does not
extend to cover the rational connectedness property, which  we will consider in this paper.
Let us recall  the  definition (cf \cite{campana}, \cite{kmm}, \cite{kollarbook}).
\begin{Defi} A compact K\"ahler manifold $X$ is rationally connected if
for any two points
$x,\,y\in X$, there exists a  (maybe singular) rational curve
$C\subset X$ with the property that $x\in C,\,y\in C$.
\end{Defi}
Examples of rationally connected varieties are given by smooth Fano varieties, i.e.
smooth projective varieties $X$ satisfying the condition that
$-K_X$ is ample. (This is the main result of \cite{campana}, and \cite{kmm}.)

The following conjecture appears in \cite{kollarvrai}. It was
 asked to me  by Pandharipande and Starr :

\begin{conj}\label{qq} (Koll\'ar) Assume $X$ is rationally connected. Let $Y$ be  a compact K\"ahler
manifold  symplectically equivalent to $X$. Then  $Y$ is also rationally connected.
\end{conj}
\begin{rema}{\rm A compact K\"ahler manifold $X$ whis is rationally connected satisfies
$H^2(X,\mathcal{O}_X)=0$, hence is projective. Thus, under
the assumption above, $X$ is projective, and if the answer to conjecture \ref{qq} is positive,
$Y$ is also projective.}
\end{rema}
This conjecture has an easy positive answer in the case of surfaces, as an immediate consequence
of theorem \ref{thkr}. Indeed, let $X$ be rationally connected of dimension $2$, and let
$Y$ be symplectically equivalent to $X$. Then $Y$ is uniruled, as $X$ is.
On the other hand $b_1(Y)=0$, because $b_1(X)=0$ and
$Y$ is diffeomorphic to $X$. Thus $Y$ is a rational surface, hence rationally connected.

In this note, we prove  the following partial results  concerning conjecture \ref{qq} in dimension $3$.
I should mention here that in these form the results are partly due to Jason Starr.
 Indeed, in the original version
of this paper, I had worked with a more restricted notion of symplectic equivalence between compact K\"ahler
manifolds, where I considered only symplectic diffeomeorphisms
$(X,\alpha)\cong (Y,\beta)$ where $\alpha$ {\it and} $\beta$ were taming for the complex structure.
Jason Starr showed me how to make the proof of proposition
\ref{propintro} work as well  when only $\alpha$ is taming, and $\beta$ is any
symplectic form which is a deformation (through a family of symplectic forms) of a
K\"ahler form on $Y$.

\begin{prop}\label{propintro} Let $X$ be  rationally connected of dimension
$3$, and let $Y$ be  compact K\"ahler symplectically equivalent to
$X$. If  $Y$ is not rationally
connected,  $X$ and $Y$ admit  almost holomorphic rational maps
$$\phi:X\dashrightarrow \Sigma,\,\phi':Y\dashrightarrow \Sigma'$$
to a surface, with rational fibers $C$, resp. $D$, of the same class (where we use
the symplectomorphism
$\psi:X\cong Y$ giving symplectic equivalence to
identify $H_2(X,\mathbb{Z})$ and $H_2(Y,\mathbb{Z})$).
\end{prop}

Here {\it almost holomorphic} means that the map is well-defined near a generic fiber.
We then consider the case where the above map $\phi$ is well-defined.
\begin{prop} \label{casouphiestdefinie} Under the same assumptions as
in proposition \ref{propintro}, assume that the rational map
$\phi$ above is well-defined and that either $\Sigma$ is smooth, or
$\phi$ does not contract a divisor to a point. Then $Y$ is also rationally connected.
\end{prop}
We will use this result together with some birational geometry arguments to prove the following:
\begin{theo}\label{thmain} Let $X, \,Y$ be  compact K\"ahler  $3$-folds. Assume that
$X$ and $Y$ are symplectically equivalent and that
one of the two following assumptions hold:
\begin{enumerate}
\item \label{16dec}
$X$ is Fano.
\item\label{26dec} $X$ is rationally connected, and  $b_2(X)\leq 2$.
\end{enumerate}
Then $Y$ is rationally connected.
\end{theo}
This anwers conjecture \ref{qq} when $X$ is a Fano threefold or
satisfies $b_2\leq2$.
The two considered cases have a small overlap. In the class where
$b_2(X)\leq 2$, one has all the blow-ups of Fano manifolds with $b_2=1$ along a connected  submanifold.
Thus this is not a bounded family. It is known on the contrary that Fano manifolds
form a bounded family (see \cite{campana}, \cite{kmm}, or \cite{classmorimukai} for the $3$-dimensional case).
However  the bound for $b_2$ of a Fano threefold is $10$ (cf \cite{classmorimukai}), showing that the Fano case
is far from being included in the second case.
\begin{rema}{\rm  Note that for varieties with $b_2=1$, conjecture \ref{qq}
obviously has an affirmative answer. Indeed a uniruled
projective manifold with $b_2=1$ is necessarily Fano. Hence if
$X$ is rationally connected with $b_2=1$,  by theorem \ref{thkr}   any
projective manifold which is symplectomorphic to it is also uniruled with $b_2=1$, hence Fano, hence
rationally connected.}
\end{rema}

To conclude this introduction, let us sketch the proof of theorem \ref{thkr}, and explain on an example
the difficulty one meets to extend it to the rational connectedness question.

\vspace{0.5cm}

{\bf Proof of theorem \ref{thkr}.}
Let $\alpha$ be a taming symplectic form on $X$ (one can take here
a  K\"ahler form). We will denote in the sequel
the degree of curves $C$  in $X$  with respect to $\alpha$ (that is the integrals $\int_C\alpha$)
 by $deg_\alpha(C)$.
Let $\mu_\alpha(X)$ be the minimum of the following set:
$$S_X:=\{deg_\alpha(C),\,C\,\,{\rm moving\,\,rational\,\,curve\,\,in}\,\,X\}.$$
Here by ``moving'', we mean that the deformations of $C$ sweep-out $X$. Note that
the minimum of the set $S_X$ is well defined, because there are finitely many families of curves
of bounded degree in $X$ and the $(1,1)$-part $\alpha^{1,1}$ of $\alpha$
is $>\epsilon\omega$ where $\omega$ is any K\"ahler form on $X$.
Let now $C$ be a moving rational curve on $X$, which
satisfies $deg_\alpha(C)=\mu_\alpha(X)$ and let $[C]\in H_2(X,\mathbb{Z})$ be its homology class.
We claim that for $x\in X$, and for adequate cohomology classes
$A_1,\ldots,A_r\in H^4(X,\mathbb{Z})$, the Gromov-Witten invariant
$GW_{0,[C]}([x],A_1,\ldots,A_r)$ counting genus $0$ curves passing through $x$ and
meeting representatives $B_i$ of the homology classes Poincar\'{e} dual to $A_i$, is non zero.
To see this, we observe that by minimality of $deg_\alpha(C)$, any genus $0$ curve
of degree $<deg_\alpha(C)$ is not moving, that is, its deformations do not sweep-out
$X$. It follows that for a general point $x\in X$, any
genus $0$ curve
of class $[C]$ and passing through $x$ is irreducible, with normal bundle generated by sections.
This implies that the set $Z_{x,[C]}$ of rational curves of classes $[C]$ passing through
$x$
has the expected dimension and  it is nonempty by assumption. Let $r$ be its dimension, and choose for
$A_i,\,1\leq i\leq r$ a class $h^2$, where $h$ is ample line bundle on $X$. It is then clear that
$GW^X_{0,[C]}([x],A_1,\ldots,A_r)\not=0$, as this number is the degree of a big and nef line bundle
on $Z_{x,[C]}$.

As $Y$ is symplectically isomorphic to
$X$, (for some symplectic structures
on $X$, resp. $Y$,  in the  deformation
class determined by  K\"ahler forms,)
we conclude that $GW^Y_{0,\psi_*[C]}([y],A'_1,\ldots,A'_r)\not=0$, where
$A'_i=\psi_*A_i\in H^4(Y,\mathbb{Z})$. But in turn, because Gromov-Witten invariants can be computed
using rational curves on $Y$ by excess formulas (see \cite{liti}, \cite{befa}, \cite{siebert}), this implies that
there is through any point $y\in Y$ a rational curve of class
$\psi_*[C]$. Thus $Y$ is uniruled.
\cqfd
\begin{rema}{\rm The proof above shows in fact a strongest statement, namely the fact that
a uniruled compact K\"ahler manifold
 $X$ admits non-zero   Gromov-Witten
invariants in genus $0$ passing through one point:
$$GW^X_{0,[C]}([x],A_1,\ldots,A_r)\not=0.$$ From this point of view,
 the proof of Theorem
\ref{thmain} is somewhat different. Indeed we do not prove that a projective
rationally connected $3$-fold
  $X$ admits non-zero  Gromov-Witten in genus $0$ passing through two points:
  $GW^X_{0,[C]}([x],[x], A_1,\ldots,A_r)\not=0$, which would be the natural symplectic analogue of  rational
  connectedness.

Our argument uses Gromov-Witten invariants {\it in higher genus}, which of course
works in the symplectic setting  as well. What we show essentially is that there is
a covering family of rational curves of class $[C]$
with a non zero $1$ point Gromov-Witten invariant: $GW^X_{0,[C]}([x], A_1,\ldots,A_r)\not=0$,
and that there is a non zero Gromov-Witten invariant of the following shape
$$GW^X_{g,[C']}(\underset{r}{\underbrace{[C],\ldots,[C]}}, A_1,\ldots,A_N)\not=0,$$
for some $r>g$ and curve class $[C']$ not proportional to $C$. We have the same non vanishings for
$Y$.

The second ingredient  is  the notion
of maximal rationally connected fibration due to Koll\'ar-Miyaoka-Mori and Campana in the K\"ahler context.
This last notion does not seem to extend well to the symplectic geometry
context.
The argument consists roughly in proving that the basis of the
maximal rationally connected fibration of $Y$ cannot be a $3$-fold by the first non vanishing,
 and cannot be a surface,
which would be uniruled by the second non-vanishing. Finally it cannot be a curve by elementary
topological considerations.
}
\end{rema}
\begin{rema}{\rm We used in this sketch of proof the terminology ``rational curve in $X$'' to mean
``stable $n$-pointed genus $0$ maps $f:C\rightarrow X$'', which are the correct objects
to count in order
to compute the Gromov-Witten invariants (cf \cite{fupagromovwitten}). However, note that if $f$ is as above,
$f(C)$ is  a rational curve in the previous sense.}
\end{rema}
If we want to apply the reasoning to study rational connectedness, we are faced to the following problem:
we could as before introduce the minimal degree for which there
are rational curves in $X$ passing through any two points of $X$. On the other hand,
it might be that curves of this degree are all reducible, with one component which is highly obstructed,
so that one cannot conclude that the corresponding Gromov-Witten invariant is non zero. In fact, consider the
case of a Hirzebruch surface
$p:F\rightarrow \mathbb{P}^1$ which is a deformation (hence symplectically equivalent to)
of  a
quadric $\mathbb{P}^1\times \mathbb{P}^1$ : Let $C_0$ be a rational curve
which is a section of $p$ with sufficiently negative self-intersection : $C_0^2<-4$.
Then one has in $F$
rational curves consisting of the union of two fibers with the section $C_0$. Such curves
$C$ can be chosen so as to pass through any two points of $F$, and  we may assume
they are, among the rational curves satisfying this property,
of minimal degree with respect to an adequate polarization. On the other hand, we have $C^2<0$
and it is clear that these curves disappear under a deformation from $F$ to
$\mathbb{P}^1\times \mathbb{P}^1$. The corresponding $2$-points Gromov-Witten invariant
is $0$ in this case.

\vspace{0.5cm}

The paper is organized as follows. In section \ref{sec1}, we prove proposition
\ref{propintro}.
In section \ref{sec2}, we study the remaining case, where $X$ is an almost  conic bundle
(we mean by this that  $X$ admits a rational map $f$ to a projective surface $\Sigma$, with
generic  fiber
isomorphic to $\mathbb{P}^1$, and that the rational map $f$  is well-defined along
the generic fiber). We show that $\phi$ is actually a morphism
(for an adequate choice of birational model of $\Sigma$)
when $b_2(X)\leq2$ or $X$ is Fano, unless there are some non trivial genus $0$ Gromov-Witten invariants
of the form $GW^X_{0,[C']}([C],A_1,\ldots,A_r)$, with $[C']$ not proportional to
$[C]$. These
Gromov-Witten invariants will be used in the last section to conclude that in this  last case, $Y$ is also rationally
connected.
We also show  that when $\phi$ is well-defined,  there are many non zero Gromov-Witten invariants on $X$,
maybe not in genus $0$ however.

The proof
of theorem \ref{thmain}  uses in turn  these non zero Gromov-Witten invariants on $Y$.
It is
 completed in section \ref{sec3}.

\vspace{0,5cm}

{\bf Thanks.} It is a pleasure to acknowledge discussions with Jason Starr and Rahul Pandharipande,
which started me thinking
to this question. I thank  Dusa McDuff, Yongbin Ruan and  Johan de Jong for comments on
various versions of the paper. I am mostly indebted to Jason Starr for showing
me how to modify my original work to get   the present version of the result.
\section{\label{sec1} Study of the rationally connected fibration of $Y$}
This section has been much simplified and improved thanks to the help of Jason Starr.
In the proof of proposition \ref{propmain} below, he showed me how to
work with general symplectic equivalence, instead of restricted symplectic equivalence
as I did originally.

We will assume that $X$ is a  projective  rationally connected complex manifold, that
$Y$ is compact K\"ahler  and that $X$ and $Y$ are symplectomorphic with respect to some symplectic forms
$\alpha,\,\beta$ on $X$, $Y$ respectively, with
$\alpha$ a taming form for the complex structure on $X$
and $\beta$  in the deformation class (as a symplectic form) of a K\"ahler form on $Y$.
We will denote as before
$\psi:X\cong Y,\,\psi^*\beta=\alpha$ such a symplectomorphism.
The theory of Gromov-Witten invariants shows that
the map $\psi$ identifies the Gromov-Witten invariants of $X$ and $Y$,  computed using holomorphic curves on
$X$ and $Y$.

We start now as in the proof of Theorem \ref{thkr}. Introducing as before
moving rational curves (or rather genus $0$ stable maps)
$C$ on $X$, of minimal degree
with respect to $\alpha$, we concluded that  there
is a covering family of rational curves
(genus $0$ stable maps) in $Y$ in the class $\psi_*([C])$.

Our goal in this section is to show the following, (which implies proposition \ref{propintro}):
\begin{prop}\label{propmain} If  $Y$ is not rationally
connected, then the covering family of curves $C$ in $X$ is given by an almost holomorphic rational map
$$\phi:X\dashrightarrow \Sigma$$
to a surface, with rational fibers of class $[C]$. Furthermore, $Y$ also admits
an almost holomorphic rational map
$$\phi':Y\dashrightarrow \Sigma'$$
with rational fiber of class $[D]=\psi_*[C]$.
\end{prop}
Here almost holomorphic means that the rational map $\phi$ is well-defined along the generic fiber of $\phi$.
Equivalently, choosing a desingularization
$$\tilde{\phi}:\widetilde{X}\rightarrow\Sigma,\,\tau:\widetilde{X}\rightarrow X$$
of $\phi$, where $\tau$ is a composition of blow-ups along smooth
centers, this means that the exceptional divisors of $\tau$ do not dominate $\Sigma$.
As the fibers of this fibration are rational curves, but $X$ is not necessarily
ruled (as it may not exist a line bundle with intersection $-1$
with fibers), we will say that $X$ is an {\it almost conic bundle}.

The proof of the proposition is based on the following lemma
(here we do not distinguish the image curve and the map, as we know that the map
 is generically the normalization map):
\begin{lemm} \label{secondlemmtrain} $Y$ is rationally connected,
unless possibly if the curve $D$ above satisfies
$c_1(K_Y)\cdot[D]=-2$ and $GW_{0,[D]}([y])=1$.
\end{lemm}
{\bf Proof.} We study the maximal rationally connected fibration of $Y$, which exists
even if $Y$ is only K\"ahler by \cite{campana}, and is an almost holomorphic rational map
$$Y\dashrightarrow B.$$
Notice that $dim\,B\leq 2$ because
 $Y$ is covered by
 rational curves $D$ of class $[D]=\psi_*[C]$.
We use now the following elementary  lemma.
\begin{lemm}\label{depuisledebut} Let $X$, $Y$
be compact K\"ahler manifolds which are symplectically
equivalent. Assume $X$  is rationally connected. If
 the basis $B$ of the rationally connected fibration of $Y$ has dimension
 $\leq 1$, $Y$ is rationally
connected.
\end{lemm}
{\bf Proof.}
We know that
$H^1(X,\mathbb{C})=0$  because
$X$ is rationally connected and this obviously implies $H^0(X,\Omega_X)=0$, hence
$H^1(X,\mathbb{C})=0$ by Hodge theory. As
$Y$ is diffeomorphic to $X$,   $H^1(Y,\mathbb{C})=0$ as well.
It follows that if the basis
$B$ of the rationally connected fibration of $Y$ has dimension
 $ 1$, it  is isomorphic to $\mathbb{P}^1$.  This contradicts
  \cite{ghs}, which implies that the basis of the rationally connected fibration is not uniruled.
\cqfd

Thus we conclude that if $Y$ is not rationally connected, the basis $B$ of the
maximal rationally connected fibration of $Y$ is a surface $\Sigma'$. Furthermore
the map
$\phi':Y\dashrightarrow \Sigma'$ is almost holomorphic.
The surface $\Sigma'$ is not uniruled by \cite{ghs}, and thus any (connected) rational curve
(or rather genus $0$ map) $ f:\Gamma\rightarrow Y$ passing through a general point $y$ of
$Y$ (where we may assume, because $\phi'$ is almost holomorphic that
$\phi'$ is well-defined everywhere along the smooth connected curve
$D':=\phi^{'-1}(\phi'(y))$) must have image supported on $D'$.
It follows that $[f_*\Gamma]=m[D']$, for some $m\geq1$.

We apply this to our  covering family of rational curves $D$
(genus $0$ stable maps) in $Y$ in the class $[D]=\psi_*([C])$ and we conclude
that $\psi_*[C]=m[D']$. Next we observe that $GW_{0,[D']}^Y([y])=1$, because
the only rational curve of class $[D']$ passing through $y$ is
$D$, which is smooth with trivial normal
bundle, so that there is fact exactly one genus $0$ map $f$ of class $[D']$ passing through $y$,
and as $H^1(N_f(-y))=0$, this stable map is computed with multiplicity $1$ in
$GW_{0,[D']}^Y([y])$.

This implies that $m=1$, because  we find that
$$GW^X_{0,\frac{1}{m}[C]}([x])\not=0,$$
so that $m>1$ would contradict the minimality of $deg_\alpha(C)$.
Hence we proved that
$$[D]=[D'],\,GW_{0,[D']}^Y([y])=1.$$
Finally, as $\phi'$ is well-defined along the generic fiber $D'$, we conclude that
$N_{D'/Y}$ is trivial, which implies by adjunction that $K_{Y}\cdot D'=K_Y\cdot[D]=-2$.
Thus lemma \ref{secondlemmtrain} is proved.

\cqfd

{\bf Proof of proposition \ref{propmain}.}
Notice that, as $\psi$ is a symplectomorphism with respect to
 symplectic forms $\alpha,\,\beta$
 of $X$, resp. $Y$, which are respective deformations of
 K\"ahler forms on $X$ resp. $Y$, $\psi^*c_1(K_Y)=c_1(K_X)$. This is indeed a standard fact of symplectic
geometry: the canonical class of a symplectic manifold $X$ is an  invariant
of the deformation class of the symplectic form $\omega$ on $X$. Indeed it can be computed using
 any almost complex structure on $X$ which
is tamed by $\omega$ or a deformation of
$\omega$, the set of such almost complex
structures being connected. This almost complex structure makes the tangent bundle into a complex vector
bundle and the canonical class is minus the first Chern class of this complex vector bundle.

Furthermore, we have by assumption
$[D]=\psi_*([C])$. Thus we have
$$c_1(K_X)\cdot[C]=c_1(K_Y)\cdot[D]=-2,\,$$
$$GW^{X}_{0,[C]}([x])=GW^Y_{0,[D]}([y])=1.$$
The first equality together with the fact that the general curve passing through the point $X$ is irreducible, and thus
has globally generated normal bundle, implies that for general $x\in X$, the normal bundle of
a curve $C$ of class $[C]$ passing through $x$ is trivial, which shows that
there are finitely many such  curves through $x$, and that the set of such curves has the expected dimension $0$.
Thus the number of these curves is equal to $GW^{X}_{0,[C]}([x])$ and this is equal to $1$ by the second equality above.
In conclusion we proved that if $\Sigma_0$ is the set parameterizing rational curves
in $X$ of class $[C]$ and $\Sigma$ is the union of components of $\Sigma_0$ parameterizing moving curves, then
the universal curve $$q':\mathcal{C}\rightarrow \Sigma,\,\Phi':\mathcal{C}\rightarrow X,$$ has
the property that $\Phi'$ has degree $1$. Thus $\Phi'$ is birational, and
$$\phi:=q'\circ {\Phi'}^{-1}:X\dashrightarrow \Sigma$$
gives the desired fibration into rational curves.

In order to conclude the proof, it just remains to prove that the rational map
$\phi: X\dashrightarrow \Sigma$ is  almost holomorphic. Assume this is not the case: let
$\tau:X'\rightarrow X$ be a composition of blow-ups along smooth
centers, such that $\tilde{\phi}:=\phi\circ\tau$ is well-defined.
Assume there is an exceptional divisor $E\subset X'$ which  dominates $\Sigma$ and is contracted
to a curve $Z$ (or a point) in $X$. Then
if $C'$ is the general fiber of $\tilde{\phi}$, $C'$ meets $E$. On the other
hand,  $K_{X'}=\tau^* K_X+F$ where $F$ is an effective divisor
supported on the exceptional locus, and the multiplicity of $E$ in $F$ is $>0$.
Thus we find that
$$c_1(K_{X'})\cdot[C']=-2=(\tau^*c_1(K_X)+F)\cdot [C']$$
$$>\tau^*c_1(K_X)\cdot [C']=c_1(K_X)\cdot[C]=-2,$$
which is a contradiction.

\cqfd

\section{ \label{sec2} The case where $X$ is an almost  conic bundle}
We now study almost conic bundles $\phi: X\dashrightarrow \Sigma$ with generic
fiber $C$. When $X$ is rationally connected,
$\Sigma$ is a rational surface, and thus we may assume to begin that $\Sigma=\mathbb{P}^2$.
(Indeed, the fact that $\phi$ is almost a morphism does not depend on the birational
model of the target.) Notice that, because $\phi$ is almost holomorphic, we have
$H\cdot C=0$, where the line bundle $H$ on $X$ is defined by
$${H}:=\phi^*\mathcal{O}_{\mathbb{P}^2}(1).$$
The first result is the following:
\begin{prop} \label{interth}Assume that either $X$ is Fano, or
$b_2(X)=2$. Then $H$  is numerically effective, unless we are not in the Fano case, and
there exists a curve class $[C']$ not proportional to $[C]$ such that for some cohomology classes
$A_1,\ldots,A_r\in H^*(X)$,
$$GW^X_{0,[C']}([C],A_1,\ldots,A_r)\not=0.$$
\end{prop}
{\bf Proof.}  Suppose first that $X$ is Fano. Then any irreducible  curve
$Z\subset X$ satisfies $K_X\cdot Z<0$, hence the Chow variety of its cycle is at least  one dimensional
because $dim\,X=3$ (cf \cite{kollarbook}, theorem 1.15).
 (This can also be formulated by evaluating the dimension
 of the space of deformations of the composed map
 $\widetilde{Z}\rightarrow Z\rightarrow X$, where $\widetilde{Z}\rightarrow Z$ is the normalization.)
 Thus,  $Z$ being irreducible, its cycle can be moved so as to be
 not contained in the indeterminacy locus of $\phi$. Thus $\phi^*H\cdot Z\geq0$.

 Suppose now that $b_2(X)=2$ but $X$ is not Fano. We have to show that
 either $H$ is numerically effective, or there exists a curve class $[C']$ not proportional to $[C]$ such that for some cohomology classes
$A_1,\ldots,A_r\in H^*(X)$,
$$GW^X_{0,[C']}([C],A_1,\ldots,A_r)\not=0.$$
 As $K_X$ is not nef, there exists a Mori contraction $c:X\rightarrow X'$, with
 $(Pic\,X')\otimes \mathbb{Q}=\mathbb{Q}$ and $-K_{X/X'}$ relatively
 ample.
We consider the three possible dimensions of $X'$ (cf \cite{mori}).

1) $dim\,X'=1$, that is $X'=\mathbb{P}^1$. In this case, the contraction is given by a pencil
 whose fibers are Del Pezzo surfaces. Let $L=c^*\mathcal{O}_{\mathbb{P}^1}(1)$. If $L\cdot C=0$, then
 $L$ is proportional to $H$ (because $b_2(X)=2$), and this contradicts the fact
  that the Iitaka dimension of $H$ is at least
 $2$. In the other case, we observe that the fibers of $c$ are uniruled. Fix a polarization
 $h$ on $X$ and introduce
 the minimal degree with respect to $h$ of rational curves contained in the fibers of $c$ and sweeping-out
 $X$. Let $[C']$ be a class curve such that
 $L.[C']=0$ and  achieving this minimal degree. All
 curves of class
 $[C']$ are supported on fibers of $c$.
 Exactly as in the proof of theorem \ref{thkr}, one then shows that for a covering family
 of rational curves $C'$ of this minimal degree, the generic member is irreducible with semipositive normal
 bundle. Using now the fact that $C$ intersects non trivially  the generic fiber of $c$,
 one concludes immediately that there is a non zero Gromov-Witten invariant
 $$GW^X_{0,[C']}([C],A_1,\ldots,A_r).$$

 2) $dim\,X'=2$. We have $c^*Pic\,X'=\mathbb{Z}L$, where $L$ is ample on $X'$, and if $C\cdot L=0$ we conclude as before that
 $L$ is proportional to $H$. In this case $H$ is numerically effective. In the other case,
 the map $c:X\rightarrow X'$ has for generic fiber a rational curve $C'$ with trivial normal bundle
 and satisfying $K_{X}\cdot C'=-2$. Furthermore there are only finitely many
 $2$-dimensional fibers of $c$.
 If $C$ is generic, there is thus exactly a $1$-dimensional family of
   fibers $C'$ meeting $C$, and this is exactly the expected dimension.
 It thus follows that there is a non trivial Gromov-Witten invariant
 $$GW^X_{0,[C']}([C],A_1),$$
 where $A_i=h^2\in H^4(X,\mathbb{Z})$ for some ample class $h\in H^2(X,\mathbb{Z})$.

 3) $dim\,X'=3$. In this case $c$ is a divisorial contraction. Note that
 $C$ is not proportional to the contracted extremal ray, because $C$ is a moving curve.
  A look at the list of divisorial contractions
 (cf \cite{mipebook})
 shows the following (see \cite{ruanextremalray}): Let $E$ be the exceptional divisor of the contraction, so
 that $E$ is either a ruled surface contracted to a smooth curve, or
 $\mathbb{P}^2$ or $\mathbb{P}^1\times\mathbb{P}^1$ contracted to a point.
 Let $[C']$ be the class of the fiber of the contracting ruling in the first case, or the class
 of a line in the second case, or the class of one of the two  rulings in the third case.
 Then for any curve
 class $\gamma$ such that $\gamma\cdot E\not=0$, one has
$GW^X_{0,[C]'}(\gamma,A_1,\ldots,A_r)\not=0$, for an adequate number $r$, which will be in fact
$0$ or $1$.

On the other hand $E\cdot C=0$ is impossible, because in this case
 $E$ and $H$ would be proportional in $Pic\,X$, and $E$ is contractible while the Iitaka dimension of
 $H$ is at least $2$.
 We deduce from this that  one has
 $GW^X_{0,[C']}([C],A_1,\ldots,A_r)\not=0$,
 where $A_i=h^2\in H^4(X,\mathbb{Z})$ for some ample class $h\in H^2(X,\mathbb{Z})$.

\cqfd
From this, we get the following result:
\begin{coro}\label{corocruc}Assume that either $X$ is Fano, or
$b_2(X)=2$. Then there exists a well-defined  morphism $\phi:X\rightarrow\Sigma$ with fiber
$C$,
where $\Sigma$ is a normal surface, unless we are not in the Fano case and
there exists a curve class $[C']$ not proportional to $[C]$ such that for some cohomology classes
$A_1,\ldots,A_r\in H^*(X)$,
$$GW^X_{0,[C']}([C],A_1,\ldots,A_r)\not=0.$$

\end{coro}
{\bf Proof.} We use the contraction theorem
 (cf. \cite{kawamata}, or  \cite{mipebook}, p 162) which tells that
 such a morphism exists if and only if $H$ is numerically effective and
 the curves $Z\subset X$ satisfying
 $Z.H=0$ also satisfy $Z\cdot K_X<0$.

 Indeed, by the previous theorem, we know
that ${H}$ is numerically effective,
unless there exists a curve class $[C']$ not proportional to $[C]$ such that for some cohomology classes
$A_1,\ldots,A_r\in H^*(X)$,
$$GW^X_{0,[C']}([C],A_1,\ldots,A_r)\not=0.$$
Thus, in order to apply the contraction theorem, we just have to show that
for any curve $Z\subset X$ satisfying the condition $Z\cdot H=0$, one has
$K_X\cdot Z<0$.

In the Fano case, this is obvious.
 When $b_2(X)=2$, the orthogonal of
$H$ in $H_2(X,\mathbb{Q})$ is generated by the class of $C$, which satisfies
the condition $C\cdot K_X=-2$.
\cqfd
We will use the following observation:
\begin{lemm}\label{remarquederniere} Assuming $\phi $ is well defined and either
$b_2(X)\leq2$ or $X$ is Fano,
 we may furthermore assume (by changing $\Sigma$
if necessary) that $\phi$ does not contract a divisor to a point of $\Sigma$.
\end{lemm}
{\bf Proof.} First of all, note that if $b_2(X)\leq2$, $\phi$ cannot contract a divisor
$D$ to a point of $\Sigma$. Indeed, such a divisor would satisfy
$D\cdot C=0$, hence would be proportional to $H$. But the Iitaka dimension of $H$ is
$2$, while no multiple of $D$ moves, which is a contradiction.

Consider now the Fano case. Let $x$ be a point of $\Sigma$,
and let  $E$ be the pure $2$-dimensional
part of $\phi^{-1}(x)$, (counted with multiplicities). We claim that
$-E$ is numerically effective on the fibers of $\phi$ and non trivial on  $\phi^{-1}(x)$.

Assuming
 the claim,
$H-\epsilon E$ remains numerically effective for a sufficiently small
$\epsilon$. On the other hand, curves $Z$ satisfying $Z\cdot (H-\epsilon E)=0$ satisfy
the condition $K_X\cdot Z<0$ for the same reasons as before, hence we can apply the contraction theorem
to $H-\epsilon E$, which does not contract $E$ anymore. This leads eventually to a morphism $\phi'$ which does
not contract any divisor to a point.

To see the claim, we observe that $-E\cdot F=0$ for any irreducible
curve $F$ contained in a fiber of $\phi$ but not contained in $\phi^{-1}(x)$.
Furthermore $-E_{\mid E}$ is effective and non trivial on each component of $E$. This implies that
$-E\cdot F\geq0$ for any irreducible curve $F\subset E$ whose deformations cover a
$2$-dimensional component of $\phi^{-1}(x)$. Consider now any irreducible curve
$F\subset X$   contained in
$\phi^{-1}(x)$. As $X$ is Fano of dimension
$3$,   the cycle of any such $F$ deforms to cover at least
a divisor in $X$ (cf \cite{kollarbook}, Theorem 1.15). On the other hand,
 all such deformations remain contained in a fiber of $\phi$. It follows that
 either the cycle of $F$ deforms to cover a $2$-dimensional component of $\phi^{-1}(x)$, so that
$-E\cdot F\geq 0$ as shown previously, or the cycle of
$F$ can be moved to be supported in another fiber, in which case we have $-E\cdot F=0$.

\cqfd
We consider now the case where $\phi$ is well defined (but $\Sigma$ may be singular). Our main result
is the following:
\begin{theo}\label{thmainsec2} Let $X$ be a rationally connected $3$-fold which admits a
morphism $\phi: X\rightarrow \Sigma$ to a normal surface $\Sigma$,
with generic fiber a rational curve $C$. Assume that either $\Sigma$ is smooth, or $\phi$ does not
contract a divisor to a point of $\Sigma$.
 Then there exist integers $g$, $r$  with
$g<r$,
 cohomology classes $A_1,\ldots,A_N\in  H^4(X,\mathbb{Z})$
and a
  homology class $[C']\in H_2(X,\mathbb{Z})$ not proportional to
  $[C]$ such that
$$GW^X_{g,[C']}(\underset{r}{\underbrace{[C],\ldots,[C]}}, A_1,\ldots,A_N)\not=0.$$
\end{theo}

Before giving the proof, let us establish a few lemmas.

\begin{lemm} \label{lenoel} $\Sigma$ contains a complete linear system of generically smooth
curves $Z$  of genus $g$, which do not meet generically the singular locus
of $\Sigma$, and satisfy
\begin{eqnarray} r=h^0(\Sigma,\mathcal{O}_\Sigma(Z))-1=h^0(Z,\mathcal{O}_Z(Z))> g.\label{form27}
\end{eqnarray}
\end{lemm}

{\bf Proof.} If $\Sigma$ is smooth, $\Sigma$ is rational and the result is obvious (we can even take
$g=0$). In general,
 we start from a ``very moving'' generic  smooth rational
curve $\Gamma_0\subset X$. Recall that ``very moving''  means that the normal bundle
$N_{\Gamma_0/X}$ is ample. Using the assumption that no divisor is contracted to a point
by $\phi$ or that $\Sigma$ is smooth, one concludes that
for $\Gamma_0$  generic,
$\phi( \Gamma_0)=:\Gamma'_0$
  avoids the singular locus of $\Sigma$.

Let $\mathcal{L}:=\mathcal{O}_{\Sigma}(\Gamma'_0)$.
Observe that $H^1(\Sigma,\mathcal{O}_\Sigma)=0$, because $\Sigma$ admits a
desingularization which is rationally connected.
It follows that the restriction map:
$$H^0(\Sigma,\mathcal{L})\rightarrow H^0(\Gamma'_0,N_{\Gamma'_0/\Sigma})$$
is surjective.
Observe now that because the equisingular  deformations of $\Gamma'_0$ in $\Sigma$
(which are singular rational  curves)
cover $\Sigma$, one has
$K_{\Sigma}\cdot \Gamma'_0<0$.

In fact we may even assume
$K_{\Sigma}\cdot \Gamma'_0<-1$, replacing if necessary $\Gamma_0$ by a ramified cover of it,
which by ampleness of the normal bundle can be deformed to an embedding.

It thus follows that
$$deg\,N_{\Gamma'_0/\Sigma}=deg\,K_{\Gamma'_0}\otimes {K_\Sigma^{-1}}_{\mid \Gamma'_0}\geq deg\,K_{\Gamma'_0}+2.$$
This inequality implies that
 the linear system $H^0(\Gamma'_0,N_{\Gamma'_0/\Sigma})$ has no base-point
on $\Gamma'_0$ so that a generic deformation $Z$ of $\Gamma'_0$ is smooth.
Letting $g$ be the arithmetic genus of $\Gamma'_0$, that is the genus of a generic
deformation $Z$ of $\Gamma'_0$ in $\Sigma$,  we now find  that $Z$ satisfies the desired property
$$r=h^0(\Sigma,\mathcal{O}_\Sigma(Z))-1=h^0(Z,\mathcal{O}_Z(Z))> g=h^0(Z,K_Z),$$
because $deg\,\mathcal{O}_Z(Z)\geq deg\,K_Z+2$ by adjunction and because
$K_{\Sigma}\cdot Z<-1$.
\cqfd
\begin{rema}\label{plutard}
{\rm The inequality $deg\,\mathcal{O}_Z(Z)\geq deg\,K_Z+2$ also implies that
$h^1(Z,\mathcal{O}_Z(Z))=0$, a fact which will be used later on.}
\end{rema}
Let $x_1,\ldots, x_r$ be $r$ generic points of $\Sigma$. Then there is a unique curve
$Z\subset \Sigma$ belonging to the linear system $\mid\mathcal{L}\mid$ and passing through
$x_1,\ldots, x_r$. This curve is smooth
and by Bertini the surface $X_Z:=\phi^{-1}(Z)$ is smooth. Choose now a section $\Gamma\subset
X_Z$ of the morphism $\phi_Z:=\phi_{\mid X_Z}:X_Z\rightarrow Z$. Let $C_i:=\phi^{-1}(x_i)$.
Let us prove now the following:
\begin{lemm}\label{utile24dec}  $\mathcal{L},\,x_1,\ldots,x_r,\,\Gamma$ being as above, for any $k>0$,
any stable map $f:\Gamma_1\rightarrow X$ of class
$$[\Gamma]+k[C]$$
meeting the $r$ generic fibers $C_1,\ldots C_r$ of $\phi$
has the property that $\phi\circ f(\Gamma_1)=Z$.
\end{lemm}
{\bf Proof. } This is almost obvious. We just have to be a little careful
with the singularities  of $\Sigma$.
Let us thus introduce a desingularization $\tau:\Sigma'\rightarrow \Sigma$ of $\Sigma$. Let
$\mathcal{L}':=\tau^*\mathcal{L}$ and $\tilde{x}_1,\ldots,\tilde{x}_r$ the points of
$\Sigma'$ over the generic points $x_1,\ldots, x_r$ of $\Sigma$.

Then if $f:\Gamma_1\rightarrow X$ is a curve as above, denote by
$\widetilde{\Gamma}'_1\subset \Sigma'$ the proper transform of
$\Gamma'_1:=\phi\circ f(\Gamma_1)\subset \Sigma$ (counted with multiplicities)
in $\Sigma'$. We observe that because the class
of $f(\Gamma_1)$ is $[\Gamma]+k[C]$ and $\phi(C)$ is a point,
 $\widetilde{\Gamma}'_1$
 belongs to one of the linear systems
$$\mid \tau^*\mathcal{L}-E\mid $$
on $\Sigma'$, where $E$ is an effective divisor supported on the exceptional locus of
the desingularization.
The linear system above has dimension $\leq r$, with equality
if and only if $E$ is empty.
As $\widetilde{\Gamma}'_1$ passes through $r$ generic points
of $\Sigma'$, it follows that the linear system
$\mid \tau^*\mathcal{L}-E\mid $ has dimension $r$. Thus $E$ is empty, and
the curve $\Gamma_1'$ does not meet the singular locus of $\Sigma$.
Hence $\Gamma'_1\in\mid\mathcal{L}\mid$, and as it passes through $x_1,\ldots,x_r$, it must be equal to
$Z$.
\cqfd

 Consider the morphism $\phi_Z:X_Z\rightarrow Z$.
The smooth fiber of $\phi_Z$ is a $\mathbb{P}^1$, and the singular fibers are
chains of $\mathbb{P}^1$'s.  Note that by successive contractions of $-1$-curves not meeting $\Gamma$,
one can construct from $F_Z$ a geometrically ruled surface $X_Z^0$. The curve $\Gamma$
is then the inverse image of a curve (still denoted $\Gamma$) in $X_Z^0$.
$\Gamma$ is a section of the structural morphism
 $p:X_Z^0=\mathbb{P}(\mathcal{E})\rightarrow Z$, where
 $\mathcal{E}:=p_*\mathcal{O}_{X_Z^0}(\Gamma)$ is a rank $2$ vector bundle on
 $Z$.
 We shall denote
by $\sigma: X_Z\rightarrow X_Z^0$ such a contraction morphism. It will be convenient to choose
the following basis $E_i$ of the lattice
$$H^2(X_Z,\mathbb{Z})/\sigma^*H^2(X_Z^0,\mathbb{Z})=\sigma^*H^2(X_Z^0,\mathbb{Z})^{\perp}.$$
We factor $\sigma:X_Z\rightarrow X_Z^0$ as a sequence of $m$ blow-ups at one point.
Let $\sigma_i:X_Z\rightarrow X_Z^i$ be the successive surfaces appearing in this factorization.
Then we define for $i\geq1$, $[E_i]:=\sigma_{i}^*[E]$, where $E$ is the exceptional curve
of the blow-up $X_Z^{i}\rightarrow X_Z^{i-1}$.
The classes $[E_i]$ are effective, and they satisfy
$$[E_i]^2=-1,\,[E_i]\cdot K_{X_Z}=-1.$$

{\bf Proof of theorem \ref{thmainsec2}.} We will denote by
$j:X_Z\hookrightarrow X$ the inclusion. For a curve
$\Gamma$ contained in $X_Z$, we will denote by
$[\Gamma]_{X_Z}\in H^2(X_Z,\mathbb{Z})$ its cohomology class in $X_Z$
and $[\Gamma]\in H^4(X,\mathbb{Z})$ its cohomology class in $X$. Hence
$[\Gamma]=j_*[\Gamma]_{X_Z}$.

 Let $g, \,r,\,x_1,\ldots,\,x_r,\,Z,\,\Gamma\subset X_Z$ be as in lemmas \ref{lenoel}, \ref{utile24dec}.
Let  $C_1,\ldots, C_r$ be the generic fibers  $\phi^{-1}(x_i)$ of $\phi$. We now consider
curves (stable maps) of genus $g$ and class $[\Gamma]+k[C]$ in $X$, where
  $k$ will be chosen sufficiently large.

The expected dimension of the family of such curves is equal to
$$-K_X\cdot([\Gamma]+k[C])=2k-K_X\cdot[\Gamma]=2k+\chi(\Gamma,N_{\Gamma/X})$$
$$=2k+\chi(\Gamma,N_{\Gamma/X_Z})+\chi(Z,N_{Z/\Sigma})=
2k+\chi(\Gamma,N_{\Gamma/X_Z^0})+r$$
$$=2k+r+\chi(Z,\mathcal{E})+g-1
=2k+r+deg\,\mathcal{E}+1-g.$$

If we consider the family of such curves
meeting $C_1,\ldots, C_r$, its expected dimension is $N:=2k+deg\,\mathcal{E}+1-g$, and
by lemma \ref{utile24dec}, we know that these curves
are all  contained in a given  surface $X_Z$, where $Z$ is a generic member
of the linear system $\mid \mathcal{L}\mid$ on $\Sigma$. Note that
$N$ is the expected dimension of the space of deformations of a smooth
curve of class $[\Gamma]+k[C]$ in $X_Z$.
If  $k$  satisfies  the condition
$\Gamma^2+2k>2g$,  choose a section
$\Gamma_k$ of $X_Z\rightarrow Z$ of  class $[\Gamma]+k[C]$ in $X_Z$.

 Then as $N_{\Gamma_k/X_Z}$ has degree $>2g-2$, it satisfies
$$H^1(\Gamma_k,N_{\Gamma_k/X_Z})=0.$$
As furthermore $H^1(\Gamma_k,(N_{X_Z/X})_{\mid \Gamma_k})=H^1(N_{Z/\Sigma})=0$  by remark \ref{plutard},
one concludes that $H^1(\Gamma_k,N_{\Gamma_k/X})=0$, so that
the deformation space of $\Gamma_k$ in $X$ is locally smooth of the right
dimension $N+r$. Furthermore, if $y_1,\ldots, y_N\in \Gamma_k$ are generic,
and $D:=\{y_1,\ldots, y_N\}$,
the restriction map:
$$H^0(\Gamma_k,N_{\Gamma_k/X_Z})\rightarrow H^0(D,(N_{\Gamma_k/X_Z})_{\mid D})$$
is an isomorphism. Choosing $N$ curves
$B_1,\ldots, B_N\subset X$ meeting
$X_Z$ in $y_1,\ldots, y_N$ respectively, we find that $\Gamma_k$ is an isolated point
in the family of curves of genus $g$ meeting $C_1,\ldots, C_r$ and $B_1,\ldots,B_N$.
This gives at least one positive contribution to
$GW_{g,[\Gamma_k]}^X(\underset{r}{\underbrace{[C],\ldots,[C]}}, [B_1],\ldots, [B_N])$.

However, in order to compute the Gromov-Witten invariant  above, we need to control
all curves in $X_Z$ whose class in $X$ is equal to $[\Gamma_k]=j_*[\Gamma_k]_{X_Z}$.

From lemma \ref{utile24dec}, we know that any curve in $X$ of
class
$[\Gamma_k]$  which meets $C_1,\ldots, C_r$ is contained in $X_Z$.
In order to conclude the proof, we have to compute the contribution
to $GW_{g,[\Gamma_k]}^X(\underset{r}{\underbrace{[C],\ldots,[C]}}, [B_1],\ldots, [B_N])$
of all the families of curves $f:\Gamma_1\rightarrow X_Z$, where
$\Gamma_1$ is  (maybe nodal)  of arithmetic genus $g$, such that the class  in $X$ of
$f(\Gamma_1)$ (counted with multiplicities)
 is equal to $[\Gamma_k]$, with $k$ large.

For this, we need the following lemma
\begin{lemm}\label{kernel}  Classes in the kernel
of $j_*:H_2(X_Z,\mathbb{Z})\rightarrow H_2(X,\mathbb{Z})$ are integral combinations
of the classes
$\frac{1}{2}[C]-[E_i]$.
\end{lemm}
{\bf Proof.} $H_2(X_Z,\mathbb{Z})$ is generated  over $\mathbb{Z}$ by
 the classes  $[C]$ of the fiber
of $\phi_Z$,  the class $[\Gamma]$ of a section of $\phi_Z$ and
the classes $[E_i]$.

If $\alpha\in Ker\,j_*$, write
$$\alpha=n [C]+m[\Gamma]+\sum_in_{i}(\frac{1}{2}[C]-[E_i]),\,n,\,m,\,n_{i}\in \mathbb{Z}.$$
Then we must have $m=0$ because $\phi_*(j_*\alpha)=0=m[Z]$.
Next we have $K_X\cdot [E_i]=-1$, because $K_{X_Z}\cdot [E_i]=-1$ and
$K_X$ has the same restriction as $K_{X_Z}$ on the fibres of $\phi_Z$.
  Furthermore $K_X\cdot [C]=-2\not=0$, and $K_X\cdot(\frac{1}{2}[C]-[E_i])=0$.
Thus
$$j_*\alpha=0\Rightarrow K_X\cdot\alpha=0\Rightarrow n=0.$$
Hence we proved that $\alpha$ is a combination of the $\frac{1}{2}[C]-[E_i]$ with integral coefficients
$n_i$.
Note that if such a combination belongs to $H_2(X_Z,\mathbb{Z})$, the
 $n_i\in \mathbb{Z}$  satisfy the condition that $2$ divides
 $\sum_in_i$.

\cqfd

We need thus to study  maps
$f:\Gamma_1\rightarrow X_Z$
where $\Gamma_1$ is a nodal curve of genus $g$,
$f_*[\Gamma_1]_{fund}=\gamma:=[\Gamma_k]+\sum_in_i(\frac{1}{2}[C]-[E_i])$.
 Note that for each such map,
  $\phi_Z\circ f:\Gamma_1\rightarrow Z$
is an isomorphism on
the (unique) genus $g$ component
of $\Gamma_1$ and contracts
all the other  components of $\Gamma_1$, which must be rational. As $deg\, N_{Z/\Sigma}>2g-2$, it follows
that $H^1(\Gamma_1, f^*N_{Z/\Sigma})=0$,
and as an easy consequence,  for fixed $\gamma$, the contribution
 of this family to
$GW_{g,[\Gamma_k]}^X(\underset{r}{\underbrace{[C],\ldots,[C]}}, [B_1],\ldots, [B_N])$
is equal to
$$GW_{g,\gamma}^{X_Z}( [B_1]_{\mid{X_Z}},\ldots, [B_N]_{\mid X_Z}).$$

Of course $[B_i]_{\mid{X_Z}}$ is a multiple of the class of a point of $X_Z$.
It thus remains to prove that for $k$ large enough
and  any $\gamma=[\Gamma_k]+\sum_in_i(\frac{1}{2}[C]-[E_i])$,
$$GW_{g,\gamma}^{X_Z}( \underset{N}{\underbrace{[pt],\ldots, [pt]}})\geq0.$$
Note that by deforming $X_Z$, we may assume the successive blow-ups starting from
$X_Z^0$ are at $m$ distinct points $z_1,\ldots, z_m\in X_Z^0=\mathbb{P}(\mathcal{E}) $.

We have the following:
\begin{lemm} \label{encoreunlemme} $m$ being fixed, for $k$ sufficiently large, for a
fixed  choice of distinct points
$z_1,\ldots, z_m\in \mathbb{P}(\mathcal{E})$, for any
choice of integers $n_1,\ldots,n_m\in\mathbb{Z}$, any linear system $L$
on the surface $X_Z'$ which is $\mathbb{P}(\mathcal{E})$ blown-up at $z_1,\ldots, z_m$, of class
$$c_1(L)=\gamma=[\Gamma_k]+l[C]-\sum_in_i[E_i],\,l=\frac{1}{2}\sum_in_i$$
satisfies $h^0(X_Z',L)\leq N+1-g$, and when equality holds, the generic member of this linear system
is smooth.
\end{lemm}
Assuming this lemma, it follows that for each $\gamma$ as above,
the dimension of the space of divisors in $X_Z'$ of class $\gamma$ has dimension $\leq N$.
Thus the dimension of the space of divisors of class
$\gamma$ passing through $N$ generically chosen points is
$0$. Furthermore, when equality holds, the finitely many divisors
of class
$\gamma$ passing through $N$ generically chosen points are smooth. It follows that
the stable maps $f:\Gamma_1\rightarrow X_Z'$  of class $f_*[\Gamma_1]_{fund}=\gamma$
passing through $N$ generically chosen points
have finitely many possible images which are
 smooth curves of genus $g$. Thus each of these $f$'s must be an isomorphism, and there are
also finitely many such stable maps $f$. It follows that
$GW_{g,\gamma}^{X'_Z}( \underset{N}{\underbrace{[pt],\ldots, [pt]}})\geq0$.
The proof of Theorem \ref{thmainsec2} is thus finished, modulo the proof of lemma
\ref{encoreunlemme}.
\cqfd
{\bf Proof of lemma \ref{encoreunlemme}.} Note that if $n_i\leq 0$, $n_iE_i$ is contained in the
fixed part of $\mid L\mid $. Thus it suffices to prove the result assuming
$n_i\geq 0$, and $l\leq\frac{1}{2}\sum_in_i$.
Next, note that because $\gamma\cdot[C]=1$, any section of $L$ vanishing to order $n_i$ at $z_i$
vanishes to order $n_i-1$ along the fiber $C_{z_i}$ passing through $z_i$.
This way, we are now reduced to the case where $n_i=0$ or $n_i=1$, and
$l\leq \frac{1}{2}\sum_in_i$.
Notice that, in both reduction steps, if either one of the $n_i<0$ or $n_i\geq2$, the inequality
becomes a   strict inequality.

We have thus  to show that for $k$
large enough, for any choice of  $s$ points  $z_{i_1},\ldots, z_{i_s}$ among $z_1,\ldots,z_m$,
for $L\in Pic\,X'_Z$, with
$$c_1(L)=[\Gamma]+k[C]+l[C]-\sum_{j\leq s}[E_{i_j}],\,l\leq \frac{s}{2},$$  we have
$h^0(X'_Z, L)\leq N+1-g$, while for $l< \frac{s}{2}$, we have $h^0(X'_Z, L)< N+1-g$.
Note that for $l=0,\,s=0$, we can take for $L$ the line bundle
$\mathcal{O}_{X_Z}(\Gamma_k)$ which has $N+1-g$ sections.

The points $z_{i_j}\in \mathbb{P}(\mathcal{E})$ determine a vector
bundle $\mathcal{E}'$ on $Z$, defined as the kernel of
the evaluation map $p_*\mathcal{O}_{\mathbb{P}(\mathcal{E})}(1)=\mathcal{E}\rightarrow \oplus
\mathcal{O}(1)_{\mid z_{i_j}}$.  Then sections of
$L$ on $X'_Z$ identify via $p_*$ to sections of $\mathcal{E}'(D)$ on $Z$, for some $D\in Pic^{k+l}(Z)$.
There are finitely many bundles $\mathcal{E}'$, and thus for $k$ large enough, and any
$l\geq0,\,deg\,D=k+l$, we have $H^1(Z,\mathcal{E}'(D))=0$.
As $deg\,\mathcal{E}'=deg\,\mathcal{E}-s$, it follows that
$$h^0(Z,\mathcal{E}'(D))=\chi(Z,\mathcal{E}'(D))=deg\,\mathcal{E}'(D)+2-2g$$
$$=deg\,\mathcal{E}-s+2k+2l+2-2g\leq deg\,\mathcal{E}+2-2g+2k=h^0(X_Z,\Gamma_k)=N+1-g,$$
with equality only when $2l=s$.

When equality holds, we have seen
that all the $n_i$ must be equal to $0$ or $1$, and the fact that the generic curve of class
$\gamma$ is smooth is deduced from the fact that with the notation above, the bundle
$\mathcal{E}'(D)$ is generated by sections, for $D\in Pic^{k+l}(Z)$.

\cqfd
\section{\label{sec3} Proofs of the main results}
{\bf Proof of Proposition \ref{casouphiestdefinie}.}
Here $\psi:X\cong Y$ is a  symplectomorphism with respect to some  symplectic forms $\alpha,\,\beta$ on
$X$, resp. $Y$, where $\alpha$ tames the complex structure on $X$
and $\beta$ is a deformation (as a symplectic form) of a K\"ahler form
on $Y$.
 We assume that the conclusion of proposition \ref{propintro} holds, but furthermore
the rational map $\phi:X\dashrightarrow\Sigma$ is well-defined, and that either
$\phi$ does not contract a divisor, or $\Sigma$ is smooth.
We can thus apply the conclusion of Theorem \ref{thmainsec2}. This tells us
 that
there exist integers $g<r$,
 cohomology classes $A_1,\ldots,A_N\in  H^4(X,\mathbb{Z})$
and a
  homology class $[C']\in H_2(X,\mathbb{Z})$ not proportional
  to $[C]$ such that
$$GW^{X}_{g,[C']}(\underset{r}{\underbrace{[C],\ldots,[C]}}, A_1,\ldots,A_N)\not=0.$$
It follows that the curve class $[D']=\psi_*[C']$ and the cohomology classes
$A_i:=\psi_*A_i\in H^*(Y)$ satisfy:
$$
GW^{Y}_{g,[D']}(\underset{r}{\underbrace{[D],\ldots,[D]}}, A'_1,\ldots,A'_N)\not=0.$$
But then this means that there exist a curve $D'$ of genus $g$  in $Y$, of class
not proportional to $[D]$, meeting $r$ generic
  fibers $D_1,\ldots, D_r$ of $\phi'$. This implies that  the surface
  $\Sigma'$ contains genus $g$ curves $D'':=\phi'(D')$ passing through
  $r$ generic points, with $r>g$. In fact we will rather consider in the following lemma
  these curves as stable maps
  from a nodal curve to $\Sigma$. The normal bundle should be thought as $N_{\phi'}$.
  \begin{lemm}\label{new24} If $\Sigma'$ satisfies this property, the Kodaira dimension of
  $\Sigma'$ is $-\infty$.
  \end{lemm}
  {\bf Proof.} Indeed, the generic curve $D''$ above has genus $g$ and
   satisfies
  $$h^0(N_{D''/\Sigma'}/Tors)\geq r>g,$$ where $Tors$ is the torsion of
  $N_{D''/\Sigma'}$.
  It follows that
  $D''$ contains at least one moving irreducible component
  $D''_0$ which  has
  genus $g_0$, and satisfies $$h^0(D''_0,N_{D''_0/\Sigma'}/Tors)>g_0.$$
  We claim that this implies $deg\,(N_{D''_0/\Sigma'}/Tors)>2g_0-2$. Assuming the claim, it follows
  that $deg\,(N_{D''_0/\Sigma'})>2g_0-2$, hence by adjunction that
  $K_{\Sigma'}\cdot D''_0<0$.
This implies that $h^0(\Sigma',{K_{\Sigma'}^{\otimes l}}_{\mid D''_0})=0,\forall l>0$,
  and as $D''_0$ is moving, this implies that $h^0(\Sigma',{K_{\Sigma'}^{\otimes l}})=0,\forall l>0$.

To see the claim, observe that Riemann-Roch gives
$$ \chi(D''_0,N_{D''_0/\Sigma'}/Tors)=deg\,(N_{D''_0/\Sigma'}/Tors)+1-g_0.$$
Thus, if $deg\,(N_{D''_0/\Sigma'}/Tors)\leq 2g_0-2$ and
$h^0(D''_0,N_{D''_0/\Sigma'}/Tors)>g_0$, we find that $h^1(D''_0,N_{D''_0/\Sigma'}/Tors)\not=0$.
Thus by Serre duality, $$h^0(D''_0,(N_{D''_0/\Sigma'}/Tors)^*\otimes K_{D''_0})\not=0.$$
But then this implies, because $D''_0$ is irreducible, that
$$h^0(D''_0,N_{D''_0/\Sigma'}/Tors)\leq h^0(D''_0,K_{D''_0})=g_0,$$
which is a contradiction.
\cqfd
Thus we conclude in this case that $\Sigma'$ is (birationally) a ruled surface, and it follows
that the basis of the rationally connected
  fibration
 of $Y$ has dimension $\leq1$. By  lemma \ref{depuisledebut}, $Y$ is rationally connected.
\cqfd
{\bf Proof of theorem \ref{thmain}.} We assume that $X$ and $Y$ are symplectically equivalent
and that, either $X$ is Fano, or
$X$ is rationally connected with $b_2(X)\leq2$. Thus
there is a  symplectomorphism $\psi: X\cong Y$
between $X$ endowed with a K\"ahler form $\alpha$ and
$Y$ endowed with a symplectic form  $\beta$ which is a deformation of  K\"ahler form.

We want to show that $Y$ is rationally connected.
We argue by contradiction, and assume
 that $Y$ is not rationally connected. Applying lemma \ref{secondlemmtrain}, we find that
there are  curve classes $[C),\,[D]$ on $X$ resp. $Y$, satisfying the following properties:
\begin{enumerate}
\item \label{1lundi10}$c_1(K_{Y})\cdot [D]=-2=c_1(K_X)\cdot[C]$.
\item \label{2lundi10}$GW^{Y}_{0,[D]}([y])=1=GW^{X}_{0,[C]}([x])$.
\item \label{3lundi10} The class $[C]$ is of minimal degree
with respect to $\alpha$, among those class curves satisfying
the property $GW^{X}_{0,[C]}([x])\not=0$.
\end{enumerate}
Furthermore, as proved in proposition \ref{propmain}, the manifolds $X$ and  $Y$
 are in this case
 almost conic bundles
with fiber $D$, resp. $C$ of class $[D]$, resp. $[C]$ where $[D]=\psi_*[C]$. Let us denote
by $\phi: X\dashrightarrow \Sigma$,
and $\phi':Y\dashrightarrow \Sigma'$ the almost conic bundle structures on $X$ and $Y$
respectively.

 Our assumption is that $b_2(X)\leq2$ or $X$ is Fano.
Hence we can apply to $X$ the corollary \ref{corocruc}, because $X$ is an almost conic bundle
with fiber $C$. Thus we conclude, with the notations of this section, that
 the morphism $\phi:X\dashrightarrow \Sigma$ with fiber $C$
 is well-defined, unless
there exists a curve class $[C']$ not proportional to $[C]$ such that for some cohomology classes
$A_1,\ldots,A_r\in H^*(X'_{n-1})$,
$$GW^{X}_{0,[C']}([C],A_1,\ldots,A_r)\not=0.$$
However, in the later  case, we conclude, by denoting $[D']=\psi_*[C']$, $A'_i=\psi_*A_i$, that
$$GW^{Y}_{0,[D']}([D],A'_1,\ldots,A'_r)\not=0.$$
 It follows  that there exists a rational curve of class $[D']$ which meets
 a generic curve $D\subset Y$ and as  $[D']$ is not proportional to
 $D$, we conclude that $\phi'(D')$ is not a point. Hence it follows that
 the surface $\Sigma'$ is swept-out by rational curves and the basis of the rationally connected
  fibration
 of $Y$ has dimension $\leq1$, which implies by
 lemma \ref{depuisledebut} that $Y$ is rationally connected, a contradiction.

 Thus  the morphism $\phi:X\rightarrow \Sigma$ with fiber $C$
 is well-defined. Furthermore, by lemma
 \ref{remarquederniere}, we may assume that
 $\phi$ does not contract a divisor to a point.
 By proposition \ref{casouphiestdefinie},  $Y$ is then  rationally connected,
 which is a contradiction.


\begin{thebibliography}{10}

\bibitem{befa}
Kai Behrend and Barbara Fantechi.
\newblock The intrinsic normal cone.
\newblock {\em Invent. Math.}, 128(1):45--88, 1997.

\bibitem{campana}
Fr{\'e}d{\'e}ric Campana.
\newblock Connexit{\'e} rationnelle des vari{\'e}t{\'e}s de {F}ano.
\newblock {\em Ann. Sci. {\'E}cole Norm. Sup. (4)}, 25(5):539--545, 1992.

\bibitem{cape} Fr{\'e}d{\'e}ric Campana, Thomas Peternell.
\newblock Rigidity theorems for primitive Fano threefolds.
\newblock {\em Comm. Analysis Geom.} 2 (1994), 173-201.

\bibitem{fupagromovwitten}
William Fulton and Rahul Pandharipande.
\newblock Notes on stable maps and quantum cohomology.
\newblock In {\em Algebraic geometry---Santa Cruz 1995}, volume~62 of {\em
  Proc. Sympos. Pure Math.}, pages 45--96. Amer. Math. Soc., Providence, RI,
  1997.

\bibitem{ghs}
Tom Graber, Joe Harris, and Jason Starr.
\newblock Families of rationally connected varieties.
\newblock {\em J. Amer. Math. Soc.}, 16(1):57--67 (electronic), 2003.



\bibitem{ruan2}
Jianxun Hu, Tian-Jun Li, and Yongbin Ruan.
\newblock Birational cobordism invariance of uniruled symplectic manifolds.
\newblock {\em Invent. Math.}, 2008.

\bibitem{kawamata} Yujiro Kawamata, Katsumi Matsuda, Kenji Matsuki.
\newblock Introduction to the minimal model problem.
\newblock {\em Algebraic geometry, Sendai 1985}, 283-360, Adv. Stud. Pure Math. 10, North-Holland, Amsterdam (1987).

\bibitem{kollarvrai} J{\'a}nos Koll{\'a}r.
\newblock Low degree polynomial equations: arithmetic, geometry and topology.
\newblock {\em  European Congress of Mathematics, Vol. I}, (Budapest, 1996), 255--288,
Progr. Math., 168, Birkh{\"a}user, Basel, 1998.

\bibitem{kollarbook}
J{\'a}nos Koll{\'a}r.
\newblock {\em Rational curves on algebraic varieties}, volume~32 of {\em
  Ergebnisse der Mathematik und ihrer Grenzgebiete. 3. Folge. A Series of
  Modern Surveys in Mathematics [Results in Mathematics and Related Areas. 3rd
  Series. A Series of Modern Surveys in Mathematics]}.
\newblock Springer-Verlag, Berlin, 1996.

\bibitem{kmm}
J{\'a}nos Koll{\'a}r, Yoichi Miyaoka, and Shigefumi Mori.
\newblock Rational connectedness and boundedness of {F}ano manifolds.
\newblock {\em J. Differential Geom.}, 36(3):765--779, 1992.

\bibitem{komimo2}
J{\'a}nos Koll{\'a}r, Yoichi Miyaoka, and Shigefumi Mori.
\newblock Rationally connected varieties.
\newblock {\em J. Algebraic Geom.}, 1(3):429--448, 1992.



\bibitem{liti}
Jun Li and Gang Tian.
\newblock Comparison of algebraic and symplectic {G}romov-{W}itten invariants.
\newblock {\em Asian J. Math.}, 3(3):689--728, 1999.

\bibitem{mcduff1}
Dusa McDuff.
\newblock The structure of rational and ruled symplectic {$4$}-manifolds.
\newblock {\em J. Amer. Math. Soc.}, 3(3):679--712, 1990.

\bibitem{mcduff2}
Dusa McDuff.
\newblock Hamiltonian {$S^1$} manifolds are uniruled, preprint.
\newblock 2007.

\bibitem{mipebook}
Yoichi Miyaoka and Thomas Peternell.
\newblock {\em Geometry of higher-dimensional algebraic varieties}, volume~26
  of {\em DMV Seminar}.
\newblock Birkh\"auser Verlag, Basel, 1997.

\bibitem{mori}
Shigefumi Mori.
\newblock Threefolds whose canonical bundles are not numerically effective.
\newblock {\em Ann. of Math. (2)}, 116(1):133--176, 1982.

\bibitem{classmorimukai}
Shigefumi Mori and Shigeru Mukai.
\newblock Classification of {F}ano {$3$}-folds with {$B\sb{2}\geq 2$}.
\newblock {\em Manuscripta Math.}, 36(2):147--162, 1981/82.

\bibitem{ruanextremalray}
Yongbin Ruan.
\newblock Symplectic topology and extremal rays.
\newblock {\em Geom. Funct. Anal.}, 3(4):395--430, 1993.

\bibitem{ruan1}
Yongbin Ruan.
\newblock Virtual neighborhoods and pseudo-holomorphic curves.
\newblock In {\em Proceedings of 6th G\"okova Geometry-Topology Conference},
  volume~23, pages 161--231, 1999.
\bibitem{siebert} Bernd Siebert.
\newblock Algebraic and symplectic Gromov-Witten invariants coincide,
 \newblock {\em Ann. Inst. Fourier} 49 (1999) 1743--1795.
\end{thebibliography}
\end{document}